\documentclass{amsart}

\usepackage{amsmath,amssymb,amsthm,tabularx,rotating,setspace}
\usepackage{diagrams}
\usepackage{hyperref}

\input cyracc.def 
\font\tencyr=wncyr10
\def\russe{\tencyr\cyracc}
\def\Sha{\text{\russe{Sh}}}

\newcommand{\QQ}{\mathbb Q}
\newcommand{\FF}{\mathbb F}
\newcommand{\ZZ}{\mathbb Z}
\newcommand{\BB}{\mathbb B}
\newcommand{\NN}{\mathbb N}

\DeclareMathOperator{\Sel}{Sel}
\DeclareMathOperator{\coker}{coker}
\DeclareMathOperator{\Hom}{Hom}

\DeclareMathOperator{\GL}{GL}
\DeclareMathOperator{\cork}{cork}
\DeclareMathOperator{\Gal}{Gal}

\DeclareMathOperator{\cd}{cd}
\DeclareMathOperator{\rank}{rank}
\DeclareMathOperator{\Aut}{Aut}
\DeclareMathOperator{\End}{End}
\DeclareMathOperator{\ad}{ad}
\DeclareMathOperator{\Ad}{Ad}
\DeclareMathOperator{\Lie}{Lie}

\newcommand{\Sc}{\Sel(E/F^{cyc})}

\newcommand{\tE}{\tilde{E}_{v,p^\infty}}
\newcommand{\Ep}{E_{p^\infty}}
\renewcommand{\SE}{\Sel(E/F_\infty)}

\newcommand{\Fiw}{F_{\infty, w}}
\newcommand{\kiw}{k_{\infty,w}}

\newcommand{\Ap}{A_{p^\infty}}
\newcommand{\Fc}{F^{cyc}}

\newcommand{\Fcw}{F^{cyc}_w}
\newcommand{\Fcv}{F^{cyc}_w}

\newcommand{\tEp}{\tilde{E}_{v,p^\infty}}

\newcommand{\SEc}{\Sel'(E/F^{cyc})}
\newcommand{\fH}{\mathfrak{H}}

\newcommand{\fJ}{\mathfrak{J}}

\newcommand{\VE}{V_p(E)}
\newcommand{\fM}{\mathfrak{M}}
\newcommand{\fD}{\mathfrak{D}}
\renewcommand{\Im}{\text{Im}}

\newcommand{\bQp}{\bar{\QQ}_p}
\newcommand{\mup}{\mu_{p^\infty}}

\newcommand{\Qp}{\QQ_p}
\newcommand{\Zp}{\ZZ_p}
\newcommand{\QZ}{\Qp/\Zp}

\newcommand{\calC}{{\mathcal C}}

\newcommand{\MHG}{{\mathfrak M}_H(\Sigma)}
\newcommand{\vs}{\vspace{1ex}}

\DeclareMathOperator{\cyc}{cyc}

\newtheorem{thm}{Theorem}[section]
\newtheorem{prop}[thm]{Proposition}
\newtheorem{lem}[thm]{Lemma}
\newtheorem{cor}[thm]{Corollary}

\theoremstyle{definition}
\newtheorem{definition}[thm]{Definition}

\newtheorem*{remark*}{Remark}
\newtheorem*{claim*}{Claim}

\newarrow{Mapsto}|---> 

\begin{document}  

\title{Generalised Euler characteristics of Selmer groups}

\author{Sarah Livia Zerbes}

\address{Department of Mathematics\\
  Huxely Building\\
  Imperial College London\\
  London SW7 2AZ\\
}

\curraddr{Department of Mathematics\\
  Harrison Building\\
  University of Exeter\\
  Exeter EX4 4QF, UK}
\email{s.zerbes@exeter.ac.uk}

\thanks{The author is supported by EPSRC Postdoctoral Fellowship EP/F043007/1.}

\begin{abstract}
 Let $E$ be an elliptic curve defined over a number field $F$, and let $p$ be a prime $\geq 5$. 
 In this paper we study the structure of the Selmer group of $E$ over $p$-adic Lie extensions
 $F_\infty$ of $F$. In particular, under certain global and local conditions on $F_\infty$ we 
 relate the generalised $\Gal(F_\infty\slash F)$-Euler characteristic of 
 $\SE$ to the generalised Euler characteristic of the Selmer group over the cyclotomic $\ZZ_p$-
 extension of $F$. This invariant generalises the classical Euler characteristic 
 to the case when $\rank_{\ZZ}E(F)>0$.
 Moreover, we show that the global and local conditions on $F_\infty$ are satisfied for a 
 large class of $p$-adic Lie extensions of $F$.
\end{abstract}


\maketitle
\setcounter{tocdepth}{3}
\tableofcontents

\section{Introduction}\label{introduction}

 \subsection{Statement of the main results}
 
  The study of the Euler characteristics of Selmer groups of motives over $p$-adic Lie extensions is a first step towards understanding exact formulae in non-commutative Iwasawa theory. The aim of the present paper is to generalise the results in~\cite{zerbes04}, in which the motive in question is an elliptic curve of rank zero, to the case of general rank. The results are part of the author's Cambridge PhD thesis~\cite{zerbesthesis}. 
  
  Let $F$ be a finite extension of $\QQ$ and $p$ any prime number $\geq 5$. Write $F^{\cyc}$ for the cyclotomic $\ZZ_p$-extension of $F$, and put $\Gamma=\Gal(F^{\cyc}\slash F)$. Let $F_\infty$ be a Galois extension of $F$, and write $\Sigma=\Gal(F_\infty\slash F)$. We shall assume throughout that $F_\infty\slash F$  is an {\it admissible} $p$-adic Lie extension, by which we mean that it satisfies: (i) $\Sigma$ is a $p$-adic Lie group, (ii) the extension $F_\infty$ is unramified outside a finite set of primes of $F$, (iii) $F_\infty\supset F^{\cyc}$, and (iv) $\Sigma$ has no element of order $p$. Such admissible $p$-adic Lie extensions are omnipresent in in arithmetic geometry. For example, if $V$ is a finite dimensional vector space over $\QQ_p$ of dimension $n$ on which $\Gal(\overline{F}\slash F)$ acts continuously, then the splitting field $F_\infty=F_\infty(V)$ is such an admissible $p$-adic Lie extension under rather mild hypotheses, e.g. (i) $V$ is unramified outside a finite set of primes of $F$, (ii) the determinant of $F$ s some power of the cyclotomic character giving the action of $\Gal(\overline{F}\slash F)$ on the group $\mup$ of $p$-power roots of unity, and (iii) $p>n+1$ with $n\geq 2$. By a basic result of Lazard~\cite{lazard65} and Serre~\cite{serre65}, $\Sigma$ has finite $p$-cohomological dimension, which is equal to its dimension $d$ as a $p$-adic Lie group. If $X$ is any discrete $p$-primary $\Sigma$-module, we say that $X$ has finite $\Sigma$-Euler characteristic if $H^i(\Sigma,X)$ is finite for all $i$, and when this is the case, then we define 
  \[
   \chi(\Sigma,X)=\prod_{i\geq 0}\big(\# H^i(\Sigma,X)\big)^{(-1)^i}
  \]
  Now let $E$ be an elliptic curve defined over $F$, and write $\Sel(E\slash F_\infty)$ for the classical Selmer group of $E$ over $F_\infty$, endowed with its natural left action of $\Sigma$. Assume for the rest of the paper that $E$ is ordinary at $p$ in the sense that $E$ has good ordinary reduction at all places $v$ of $F$ divising $p$. Our results depend on the analogue of Mazur's conjecture for the extension $F_\infty\slash F$ which is made in~\cite{cfksv} and which we now explain. For any compact $p$-adic Lie group $G$, write $\Lambda(G)=\ZZ_p[[G]]$ for the Iwasawa algebra of $G$ (see Section~\ref{notation}). Put $H=\Gal(F_\infty\slash F^{\cyc})$, and let $\MHG$ denote the category of all finitely generated $\Lambda(\Sigma)$-modules $M$ such that $M\slash M(p)$, where $M(p)$ denotes the $p$-primary submodule of $M$, is a finitely generated $\Lambda(H)$-module. Let
  \[
   \calC(E\slash F_\infty)=\Hom\big(\Sel(E\slash F_\infty),\QQ_p\slash \ZZ_p\big),
  \]
  endowed with its natural left $\Lambda(\Sigma)$-module structure. Then it is conjectured in~\cite{cfksv}, and proven in some special cases, that we always have $\calC(E\slash F_\infty)\in\MHG$. 
  
  The essential change of hypotheses from our previous paper~\cite{zerbes04} is that we now assume that we no longer assume that the group $E(F)$ of $F$-rational points on $E$ is finite (but we do assume that the $p$-primary subgroup $\Sha(E\slash F)(p)$ of the Tate-Shafarevich group of $E$ over $F$ is finite). Thus it is no longer true that $\Sel(E\slash F_\infty)$ has finite $\Sigma$-Euler characteristic once $E(F)$ is infinite. Instead, we have to consider the generalised $\Sigma$-Euler characteristic of $\Sel(E\slash F_\infty)$ (see Section~\ref{generalizedEC2}), whose definition depends essentially on our hypothesis that $F^{\cyc}\subset F_\infty$. We can also consider the generalised $\Gamma$-Euler characteristic of $E$ over $F^{\cyc}$, whose value has been computed by the fundamental work of Perrin-Riou~\cite{perrin-riou84} and Schneider~\cite{schneider85}. In order to relate these two generalised Euler characteristics, we have to impose the following two hypotheses on the extension $F_\infty$ over $F$:

 \begin{itemize}
   \item (Fin$_{glob}$) $H^i(H,\Ep(F_\infty))$ is finite for all $i\geq 0$.
   \item (Fin$_{loc}$) For any prime $v$ of $F$ dividing $p$, $H^i(H_w,\tEp(\kiw))$ is finite for all $i\geq 0$.
 \end{itemize}

 We discuss below a general class of admissible $p$-adic Lie extensions for which (Fin$_{glob}$) and (Fin$_{loc}$) will be proven to hold. 

 \begin{thm}\label{theorem3}
  Let $F$ be a finite extension of $\QQ$, $E$ an elliptic curve defined over $F$ and $p$ a prime $\geq 5$. 
  Let $F_\infty$ be an admissible $p$-adic Lie extension of $F$ with Galois group 
  $\Sigma=\Gal(F_\infty/F)$. 
  Assume that (i) $\Sha(E/F)(p)$ is finite, (ii) $E$ is ordinary at $p$ and (iii) $\calC (E/F_\infty)\in\MHG$. 
  Also, assume that both (Fin$_{glob}$) and (Fin$_{loc}$) hold. Then $\Sel(E/F_\infty)$ has
  finite generalised $\Sigma$-Euler characteristic if and only if $\Sel(E/\Fc)$ has finite generalised
  $\Gamma$-Euler characteristic, and if this is the case, then
  \begin{equation}\label{genEC2}
   \chi\big(\Sigma,\SE\big)=\chi\big(\Gamma,\Sc\big)\times\Big|\prod_{v\in\mathfrak{M}}L_v(E,1)\Big|_p.
  \end{equation}
  Here, $\fM$ is the set of primes of $F$ not dividing $p$ whose inertia group in $\Sigma$ is infinite,
  and for a prime $v\in\fM$, $L_v(E,s)$ denotes the Euler factor of $E$ at $v$.
 \end{thm}
 
 \begin{remark*}
  If $\Sel(E/F)$ is finite, then as shown in~\cite{coatesschneidersujatha03}, $\calC (E/\Fc)$ is a $\Lambda(\Gamma)$-torsion module and $\Sc$ has finite $\Gamma$-Euler characteristic. Also, the assumption that $\calC (E/F_\infty)\in\MHG$ implies that $\calC (E/F_\infty)$ is $\Lambda(\Sigma)$-torsion (c.f.~Proposition~\ref{finitegeneration2}) and $\calC (E/\Fc)$ is $\Lambda(\Gamma)$-torsion (c.f.~Lemma~\ref{Gammatorsion}), so Theorem~\ref{theorem3} is indeed a generalisation of Theorem 1 in~\cite{zerbes04}.
 \end{remark*}

 We show in Section~\ref{examples2} below that the conditions (Fin$_{glob}$) and (Fin$_{loc}$) are satisfied for a large class of $p$-adic Lie extensions of $F$. Assume that $F_\infty$ is the fixed field of the kernel of some continuous
 representation 
 \[\rho:\Gal(\bar{F}/F)\rTo V\] 
 of $\Gal(\bar{F}/F)$ on some finite dimensional $\Qp$-vector space $V$. For a
 prime $v$ of $F$ dividing $p$, denote by $\rho_v$ the restriction of $\rho$ to $\Gal(\bar{F}_v/F_v)$. 
 Let $H=\Gal(F_\infty/\Fc)$, and denote by $\Lie(H)$ its Lie algebra.
 
 \begin{thm}\label{theorem4}
 \begin{enumerate}[(i)]
  \item If $\Lie H$ is reductive, then (Fin$_{glob}$) is satisfied. 
  \item Suppose that for every prime $v$ of $F$ dividing $p$,
  the representation $\rho_v$ is de Rham and the eigenvalues of Frobenius on 
  \[
   (V\otimes_{\Qp}\BB_{st})^{\Gal(\bar{F}_v/F_v)}   
  \]
  are Weil numbers of the same parity (c.f. Section~\ref{examples2}). Here, $\BB_{st}$ denotes Fontaine's ring of periods. Then (Fin$_{loc}$) is satisfied.
 \end{enumerate}
 \end{thm}
  
 {\bf Acknowledgements}: I am very grateful to John Coates for his support and encouragement during my PhD and for his interest in my work. I also thank David Loeffler and the referee for their careful reading of the manuscript.


 \subsection{Notation}\label{notation}
 
  For a field $F$, which will always be a finite extension of either $\QQ$ or $\QQ_l$ for 
  some finite prime $l$, we denote by $\bar{F}$ its algebraic closure. Given a finite set $S$ of primes of
  $F$, let $F_S$ be the maximal extension of $F$ contained in $\bar{F}$ which is unramified 
  outside $S$. For an extension $L$ of $F$ contained in $F_S$, we write $G_S(L)$ for the Galois group
  $\Gal(F_S/L)$. 
  
  For an odd prime $p$, denote by $\mup$ the group of $p$-power roots of unity of $\bar{F}$. Let $\Fc$
  be the cyclotomic $\Zp$-extension of $F$, which is the unique $\Zp$-extension of $F$ contained in
  $F(\mup)$. Denote the Galois group of the extension $\Fc\slash F$ by $\Gamma$.
  
  If $G$ is a compact $p$-adic Lie group, define its Iwasawa algebra $\Lambda(G)$ to be the completed group ring
  \[
   \Lambda(G)=\varprojlim \ZZ_p[G\slash U],
  \]
  where $U$ runs over all open normal subgroups of $G$.
  
  For an abelian group $M$ and $m\in\NN$, let $M_m$ be the subgroup of $M$ whose elements are
  annihilated by $m$. Given a prime $p$, define the $p$-Tate module of $M$ to be the inverse limit
  \[
   T_pM=\varprojlim M_{p^n}.
  \]
  If $M$ is $p$-primary with the discrete topology, then the Pontryagin dual  
  \[
   M^\vee=\Hom_{cts}(M,\QZ)
  \]
  of $M$ is a compact pro-$p$ abelian group.
  
  Let $G$ be a profinite group and $M$ a $G$-module. Write $M^G$ for the $G$-invariants of
  $M$. If $M$ is a discrete $p$-primary $G$-module, then the induced action of $G$ on $M^\vee$ is given by
  \[
   (\sigma f)(x)=f(\sigma^{-1}x)
  \]
  for all $\sigma\in G$. If $M$ is a discrete $G$-module, then $H^i(G,M)$ denotes the $i$th cohomology
  group formed by continuous cochains.
  
  For a number field $F$ and a finite prime $v$ of $F$, denote by $F_v$ the completion
  of $F$ at $v$. If $L$ is an algebraic extension of $F$ and $w$ a prime of $L$ above $v$, then $L_w$ is
  the union of the completions at $w$ of all finite extensions of $F$ contained in $L$. For an elliptic
  curve $E$ defined over $F$, denote by $c_v$ the Tamagawa factor of $E$ at $v$, which is defined as
  \[
   c_v=[E(F):E_0(F)],
  \]
  where $E_0(F)$ denotes the subgroup of points of $E(F)$ which have non-singular reduction. Let
  $L_v(E,s)$ be the local Euler factor of $E$ at $v$. If $E$ has good reduction at $v$, then 
  $\tilde{E}_v$ denotes the reduction of $E$ mod $v$.

\section{Preliminaries}


 \subsection{The Selmer group}\label{diagram}

  Let $S$ be a finite set of primes of $F$ containing all the primes dividing $p$, all the primes where $E$ has bad reduction and all the primes of $F$ which ramify in $F_\infty$. Note that the choice of $S$ implies that $F_\infty\subset F_S$. Let $L$ be a finite extension of $F$ contained in $F_S$. For each finite place $v$ of $F$, define   
  \[
   J_v(L)=\bigoplus_{w| v}H^1(L_w,E)(p),
  \] 
  where $w$ runs over all primes of $L$ dividing $v$. Let $\bar{L}_w$ be an algebraic closure of
  $L_w$. By choosing an embedding $\bar{L}\subset\bar{L}_w$, we can consider $\Gal(\bar{L}_w/L_w)$ as a
  subgroup of $\Gal(\bar{L}/L)$, so we get the localisation map
  \begin{equation}\label{defnlambda}
   \lambda_S(L):H^1(G_S(L),\Ep) \rightarrow \bigoplus_{v\in S}J_v(L).
  \end{equation}
  induced by restriction.
  
  For an infinite algebraic extension $K$ of $F$ contained in $F_S$, define
  \[
   J_v(K)=\varinjlim J_v(L),
  \]
  where $L$ runs over all finite extensions of $F$ contained in $K$ and the inductive 
  limit is taken with respect to the restriction maps. Also, define
  \[
   \lambda_S(K)=\varinjlim\lambda_S(L).
  \]
  The following theorem is well-known (c.f. Section 2.2 in~\cite{coates97}).
  
  \begin{thm}
   Let $K$ be an extension of $F$ contained in $F_S$. Then $\Sel(E/K)$ is given by the exact sequence
   \[
    0 \rightarrow \Sel(E/K) \rightarrow H^1(G_S(K),\Ep) \rightarrow \bigoplus_{v\in S}J_v(K).
   \]
  \end{thm}
  
  In particular, $\SE$ is given by
  \begin{equation}\label{Selmer*}
   0 \rightarrow \Sel(E/F_\infty) \rightarrow H^1(G_S(F_\infty),\Ep) \rTo^{\lambda_S(F_\infty)} \bigoplus_{v\in S}J_v(F_\infty).
  \end{equation}
  
   
 \subsection{The category $\MHG$}
 
  In~\cite{lazard65}, Lazard has proven the following important result.
  
  \begin{thm}\label{lazard}
   Let $G$ be a pro-$p$ $p$-adic Lie group with no element of finite order. Then the Iwasawa algebra $\Lambda(G)$ is left and right Noetherian and has no zero divisors.
  \end{thm}
  
  Let $R$ be an open pro-$p$ subgroup of $\Sigma$. Then Theorem~\ref{lazard} implies that $\Lambda(R)$ has a skew field of fractions $K(R)$ (c.f. Chapter 9 in~\cite{goodearlwarfield89}). 
  \vs
  
  \begin{definition} For a finitely generated $\Lambda(R)$-module $M$, define
  \[
   \rank_{\Lambda(\Sigma)}M=\frac{1}{[\Sigma:R]}\dim_{K(R)}M\otimes_{\Lambda(R)}K(R).
  \]
  Then $M$ has rank $0$ if and only if it is a torsion $\Lambda(\Sigma)$-module. 
  \end{definition}
  
  It is a simple (but crucial) observation that $\MHG$ is a subcategory of the category of
  finitely generated torsion $\Lambda(\Sigma)$-modules.
  
  \begin{prop}\label{finitegeneration2}
   Let $M\in\MHG$. Then $M$ is a torsion $\Lambda(\Sigma)$-module. 
  \end{prop}
  \begin{proof}
   It is clearly sufficient to show that $M/M(p)$ is a torsion $\Lambda(\Sigma)$-module, so we assume
   without loss of generality that $M$ is finitely generated over $\Lambda(H)$. If $M$ is not a torsion $\Lambda(\Sigma)$-module, then $\rank_{\Lambda(\Sigma)}M>0$, so $M$ contains a copy of $\Lambda(\Sigma)$. But $H$ is not of finite index in $\Sigma$, so $\Lambda(\Sigma)$ is not finitely generated over $\Lambda(H)$, giving the required contradiction.
  \end{proof}
  
  Together with the main result of Section 3.3 in~\cite{zerbes04} we obtain the following result, which
  will be important in Section~\ref{ECcalculationII2}.
  
  \begin{prop}\label{surjection}
   Suppose that $E$ has either good ordinary or split multiplicative reduction at all primes of $F$
   dividing $p$. If $\calC (E/F_\infty)\in\MHG$, then we have a short exact sequence
   \[
    0\rTo\SE\rTo H^1(G_S(F_\infty),\Ep)\rTo^{\lambda_S(F_\infty)}\bigoplus_{v\in S}J_v(F_\infty)\rTo 0.
   \]
   Here, $S$ and $J_v(F_\infty)$ are defined as in Section~\ref{diagram}.
  \end{prop} 
  
  The condition that $\calC(E/F_\infty)\in\MHG$ also has the following important consequence.
  
  \begin{lem}\label{Gammatorsion}
   If $\calC (E/F_\infty)\in\MHG$, then $\calC (E/\Fc)$ is $\Lambda(\Gamma)$-torsion.
  \end{lem}
  \begin{proof}
   If $\calC (E/F_\infty)\in\MHG$, then as shown in~\cite{cfksv}, $\SE^H$ is a cotorsion
   $\Lambda(\Gamma)$-module. Applying the snake lemma to the commutative diagram in Section~\ref{ECcalculationI2} shows that if 
   $F_\infty$ is admissible, then
   the kernel and cokernel of the natural map
   \[
    \Sc\rTo \Sel(E/F_\infty)^H
   \]
   are finitely generated $\Zp$-modules, which finishes the proof.
  \end{proof}

  
 \subsection{Generalised Euler characteristics}\label{generalizedEC2}

  Let $D$ be a discrete $p$-primary $\Gamma$-module. Then we have
  \begin{align*}
   H^0(\Gamma,D)&=D^\Gamma,\\
   H^1(\Gamma,D)&=D_\Gamma,
  \end{align*}
  and hence there is an obvious map
  \begin{align*}
   \phi_D: H^0(\Gamma,D)&\rTo H^1(\Gamma,D),\\
                     f  &\rMapsto \text{residue class of $f$}.
  \end{align*}
  For a discrete $p$-primary $\Sigma$-module $X$, define
  \[
   d_0: H^0(\Sigma,X)=H^0(\Gamma,X^H)\rTo^{\phi_{X^H}} H^1(\Gamma,X^H)\hookrightarrow H^1(\Sigma,X),
  \]
  where the last map is given by inflation. Similarly, for $i\geq 1$, define
  \[
   d_i: H^i(\Sigma,X)\twoheadrightarrow H^0(\Gamma,H^i(H,X))\rTo^{\phi_{H^i(H,X)}} H^1(\Gamma,H^i(H,X))
   \hookrightarrow H^{i+1}(\Sigma,X),
  \]
  where the first map is given by restriction. Note that this map is surjective since 
  $\cd_p(\Gamma)=1$. For convenience, define $d_{-1}$ to be the zero map. It is 
  not difficult to see that $(H^i(\Sigma,X),d_i)$ forms a complex. Denote its homology groups by $\fH_i$.
  \vspace{1ex}
  
  \begin{definition}
  The $\Sigma$-module $X$ has finite generalised $\Sigma$-Euler characteristic if 
  $\fH_i$ is finite for all $i\geq 0$. We then define 
  \[
   \chi(\Sigma,X)=\prod_{i\geq 0}\# \fH_i^{(-1)^i}.
  \]
  \end{definition}
  
  Note that this definition agrees with the usual definition of $(\Sigma,X)$ if the groups $H^i(\Sigma,X)$ are finite for all $i$. In this paper, we will only encounter the case when $H^0(\Sigma,X)$ and $H^1(\Sigma,X)$ are infinite and $H^i(\Sigma,X)$ is finite for all $i\geq 2$. In this case, $X$ has finite generalised $\Sigma$-Euler characteristic if and only if the map $d_0$ has finite kernel and cokernel, and then
  \[
   \chi(\Sigma,X)=\frac{\#\ker(d_0)}{\#\coker(d_0)}\times\prod_{i\geq 2} (\# H^i(\Sigma,X))^{(-1)^i}.
  \]
  
  To see why this is a useful generalisation of the classical $\Sigma$-Euler characteristic, 
  recall the definition of the Akashi series that first appeared in~\cite{coatesschneidersujatha03} (though not under this name) and was further studied in~\cite{cfksv}.
  
  \begin{definition} Let $M\in\MHG$.
  Then $H_i(H,M)$ is a torsion $\Lambda(\Gamma)$-module for all $i\geq 0$. 
  Let $g_{M,i}$ be its characteristic element, and define the Akashi series of $M$ to be the alternating
  product
  \[
   f_M=\prod_{i\geq 0}g_{M,i}^{(-1)^i}.
  \]
  Note that $f_M$ is only defined mod $\Lambda(\Gamma)^*$.
  \end{definition}
  
  \begin{remark*} 
   \begin{enumerate}[(i)]
    \item If $M$ is a finitely generated torsion $\Lambda(\Gamma)$-module, then its Akashi series is the same as its characteristic element.
    \item If $X$ is a discrete $p$-primary $\Sigma$-module such that $M=X^\vee\in\MHG$, then we call $f_M$ the Akashi series of $X$, and by abuse of notation, we denote it by $f_X$.
   \end{enumerate}
  \end{remark*}

  \begin{lem}\label{multiplicativity3}
   If we have a short exact sequence 
   \[
    0\rightarrow L\rightarrow M\rightarrow N\rightarrow 0
   \]
   in $\MHG$, then $f_M=f_Nf_L$.
  \end{lem}
  \begin{proof}
   See Lemma 4.1 in~\cite{coatesschneidersujatha03}.
  \end{proof}
  
  Let $M^\vee$ be the Pontryagin dual of $M$. The following proposition gives a relation between $f_M$ and the generalised $\Sigma$-Euler characteristic of $M^\vee$.

  \begin{prop}\label{akashi2}
   Let $X$ be a discrete $p$-primary $\Sigma$-module such that $X^\vee\in\MHG$. 
   If $X$ has finite generalised $\Sigma$-Euler characteristic, then the leading term of $f_X$ is $\alpha_XT^k$, where $k$ is equal to the alternating sum
   \[
    k=\sum_{i\geq 0}(-1)^i\cork_{\ZZ_p}H^i(H,X)^\Gamma
   \]
   and 
   \[
    \chi(\Sigma,X)=|\alpha_X|_p^{-1}.
   \]
  \end{prop}
  
  \begin{remark*}
   This generalizes Lemma 4.2 in~\cite{coatesschneidersujatha03}.
  \end{remark*}

  Proposition~\ref{akashi2} is a direct consequence of the following lemma.
  
  \begin{lem}\label{gammamodules2}
   Let $X$ be a finitely generated cotorsion $\Lambda(\Gamma)$-module, and let $g_X$ be a generator 
   of the characteristic ideal of $X^\vee$. If the natural map
   $\phi:X^\Gamma\rTo X_\Gamma$ has finite kernel and cokernel, then 
   $g_X(T)=T^rf(t)$, where $f(T)\in\Lambda(\Gamma)$, $T\nmid f(T)$ and
   $r=\cork_{\ZZ_p}X^\Gamma=\cork_{\ZZ_p}X_\Gamma$, and
   \[
    \frac{\#\ker(\phi)}{\#\coker(\phi)}=| f(0)|_p^{-1}.
   \]
  \end{lem}
  \begin{proof}
   To simplify notation, let $M$ denote the Pontryagin dual of $X$, and let $\psi$ be 
   the natural map $\psi:M^\Gamma \rTo M_\Gamma$. Assume
   that $\phi$ has finite kernel and cokernel. By duality, this implies that the kernel 
   and cokernel of $\psi$ are also finite, and
   $\frac{\#\ker(\phi)}{\#\coker(\phi)}=\frac{\#\coker(\psi)}{\#\ker(\psi)}$. 
   Now $M$ is $\Lambda(\Gamma)$-torsion, so using the structure
   theorem for finitely generated $\Lambda(\Gamma)$-modules, we have a pseudo-isomorphism
   \[
    M\cong\bigoplus_{i=1}^n\Lambda(\Gamma)/(g_i(T))
   \]
   for some $g_i(T)\in\Lambda(\Gamma)$, where we can use the non-canonical identification 
   of $\Lambda(\Gamma)$ with the power series ring over
   $\ZZ_p$ in one variable. We define the map $\psi$ `component-wise', so 
   we can without loss of generality assume that 
   $M=\Lambda(\Gamma)/(g(T))$. If $g(0)\neq 0$, then $M/TM$
   is finite, so as shown in Lemma 4.2 in~\cite{greenberg99}, 
   both $M^\Gamma$ and $M_\Gamma$ are finite, and the result is
   immediate. Suppose that $g(0)=0$, and write $g(T)=T^nf(T)$, where 
   $f(T)\in\Lambda(\Gamma)$ and $T\nmid f(T)$. Now $M^\Gamma$ is the kernel
   of the map $\times T:M\rTo M$, which is given by multiplication by $T$, 
   so $M^\Gamma$ is the submodule generated by $T^{n-1}f(T)+(T^nf(T))$.
   Also, we have $M_\Gamma=M/TM$, which is equal to $\Lambda(\Gamma)/(T)$ 
   since $n\geq 1$. Now $\Lambda(\Gamma)/(T)\cong\ZZ_p$, and the map
   $\psi$ is given by
   \begin{align*}
    \psi: M^\Gamma &\rTo M_\Gamma,\\
         h(T)+(T^nf(T))& \rMapsto h(0).
   \end{align*}
   It follows that if both $\ker(\psi)$ and $\coker(\psi)$ are finite, 
   then we must have $n=1$, so in fact $\ker(\psi)=0$ and
   $\#\coker(\psi)=[\ZZ_p:f(0)\ZZ_p]=| f(0)|_p^{-1}$.
  \end{proof}

  \begin{proof}[Proof of Proposition \ref{akashi2}] 
   Recall that the map $d_i$ is given by
   \begin{equation}\label{definitiondi2}
    d_i:H^i(\Sigma,A)\rTo H^i(H,A)^\Gamma \rTo^{\phi_i}H^i(H,A)_\Gamma\hookrightarrow H^{i+1}(\Sigma,A).
   \end{equation}
   It is not difficult to see that if $A$ has finite generalised $\Sigma$-Euler characteristic, then $\phi_i$ has finite kernel and cokernel for all $i\geq 0$. A careful analysis of the maps in~\eqref{definitiondi2} shows that
   \begin{align*}
    \# G_0&=\#\ker(\phi_0),\text{ and}\\
    \# G_i&=\#\ker(\phi_i)\#\coker(\phi_{i-1})
   \end{align*}
   for all $i\geq 1$, where as before $G_i=\ker(d_i)/\Im(d_{i-1})$. The proof is now immediate from Lemma~\ref{gammamodules2}.
  \end{proof}

  As shown in~\cite{rapoportzink82} the maps $d_i$ can also be constructed as follows.
  
  \begin{lem}
   Let $c$ be a generator of $H^1(\Gamma,\ZZ_p)$. The the map $d_i$ is given by
   \begin{align*}
    d_i: H^i(\Sigma,A)&\rTo H^{i+1}(\Sigma,A),\\
    f&\rMapsto f\cup c.
   \end{align*}
  \end{lem}
  
  The following result will be useful.
  
  \begin{lem}\label{sameEC2}
   Let $\phi:A\rTo B$ be a homomorphism of discrete 
   $p$-primary $\Sigma$-modules, and assume that $\phi$ has finite kernel and cokernel. Then $A$ has
   finite generalised $\Sigma$-Euler characteristic if and only if $B$ does, 
   and if this is the case, then
   \[
    \chi(\Sigma,A)=\chi(\Sigma,B).
   \]
  \end{lem}
  \begin{proof}
   We can assume without loss of generality that $\ker(\phi)=0$, so we have a short exact sequence of
   $\Sigma$-modules
   \[
    0\rTo A\rTo B\rTo C\rTo 0,
   \]
   where $C$ is finite. Taking $\Sigma$-cohomology, for all $i\geq 0$ the long exact sequence gives
   homo\-morphisms
   \[
    H^i(\Sigma,A)\rTo H^i(\Sigma,B),
   \]
   which have finite kernel and cokernel since $C$ is finite. These maps are compatible with the maps in
   the construction of the complex giving the generalised $\Sigma$-Euler characteristic, so it follows that
   $A$ has finite generalised Euler characteristic if and only if so does $B$. By the multiplicativity of
   Akashi series in short exact sequences (c.f. Lemma~\ref{multiplicativity3}), we have
   \[
    f_B=f_A\times f_C.
   \]
   But $f_C\thicksim 1$ since $C$ is finite, so the result follows immediately from
   Proposition~\ref{akashi2}.
  \end{proof}
  
  \begin{remark*}
   It is {\it not} true that the generalised $\Sigma$-Euler characteristic is multiplicative in short exact sequences. I thank Otmar Venjakob for pointing out to me the following example: identify $\Lambda(\Gamma)$ with the power series ring $\Zp[[T]]$. Then it is easy to see that in the short exact sequence
   \begin{equation}\label{counterexample2}
    0\rTo \Lambda(\Gamma)/T \rTo^{\times T} \Lambda(\Gamma)/T^2\rTo \Lambda(\Gamma)/T\rTo 0
   \end{equation}
   the first and the last non-trivial modules have finite generalised $\Gamma$-Euler characteristic, but the middle one does not. 
  \end{remark*}
  
 
   \subsection{The large Selmer group $\SEc$}\label{largeSelmer2}
  
    Let $S'$ be the subset of primes of $F$ not dividing $p$ which have infinite inertia group in $\Sigma$, 
    and let 
    $S''=S\backslash S'$. Define the large Selmer group $\SEc$ by the exact sequence
    \[
     0\rTo\SEc\rTo H^1(G_S(\Fc),\Ep)\rTo \bigoplus_{v\in S''}J_v(\Fc).
    \]
    {\bf Note.} Since $\SEc$ depends on the extension $F_\infty$, one should write $\SEc _{F_\infty}$, but we omit the index to keep the notation as simple as  possible.
    \vspace{1ex}

   It is shown in~\cite{coatesschneidersujatha03} that if ${\mathcal C}(E/\Fc)$ is $\Lambda(\Gamma)$-torsion, then $\SEc$ and $\Sel(E/\Fc)$ are 
   related by the short exact sequence
   \begin{equation}\label{largesurjective}
    0\rTo\Sc\rTo \SEc\rTo \bigoplus_{v\in S'}J_v(\Fc)\rTo 0.
   \end{equation} 
   The generalised $\Gamma$-Euler characteristics of $\Sc$ and $\SEc$ are closely related.

   \begin{lem}\label{EClarge2}
    Suppose that $\Sc$ is $\Lambda(\Gamma)$-torsion. Then
    the large Selmer group $\SEc$ has finite generalised $\Gamma$-Euler characteristic if and only if
    $\Sc$ does, and if this is the case, then
    \[
     \chi\big(\Gamma,\SEc\big)=\chi\big(\Gamma,\Sc\big)\times\Big|\prod_{v\in S'}L_v(E,1)\Big|_p.
    \]
   \end{lem}
   \begin{proof}    
    By Lemma~\ref{sameEC2} and equation~\eqref{largesurjective}, it is sufficient to 
    show that $J_v(\Fc)$ has finite
    $\Gamma$-Euler characteristic for all $v\in S'$, and that
    \[
     \chi\big(\Gamma,J_v(\Fc)\big)=\big|L_v(E,1)\big|.
    \]
    Shapiro's lemma gives canonical isomorphisms
    \begin{align*}
     H^0(\Gamma,J_v(\Fc))&\cong H^0(\Gamma_v,H^1(\Fcv,E)(p)),\\
     H^1(\Gamma,J_v(\Fc))&\cong H^1(\Gamma_v,H^1(\Fcv,E)(p)),
    \end{align*}
    where $\Gamma_v$ is the decomposition group in $\Gamma$ of some fixed prime of $\Fc$ above $v$. Since
    $v\nmid p$, Kummer theory shows that 
    \[
     H^1(\Fcv,E)(p)\cong H^1(\Fcv,\Ep).
    \]
    Applying the Hochschild-Serre spectral sequence to the extension $\Fcv$ of $F_v$ gives a short exact
    sequence
    \begin{equation}\label{H0gamma2}
     0\rTo H^1(\Gamma_v,\Ep(\Fcv))\rTo H^1(F_v,\Ep)\rTo H^0(\Gamma_v,H^1(\Fcv,\Ep))\rTo 0
    \end{equation}
    and an isomorphism
    \begin{equation}\label{H^1gamma2}
     H^2(F_v,\Ep)\cong H^1(\Gamma_v,H^1(\Fcv,\Ep)).
    \end{equation}
    Since $\Ep(F_v)$ is finite, local duality implies that $H^2(F_v,\Ep)=0$, so
    \[ 
     H^1(\Gamma_v,H^1(\Fcv,\Ep))=0.
    \]
    Also, it is shown in~\cite{coates97}, Lemmas 3.6 and 3.8, that the groups $H^1(\Gamma_v,\Ep(\Fcv))$ and 
    $H^1(F_v,\Ep)$ are finite, and that
    \begin{align*}
     \# H^1(\Gamma_v,\Ep(\Fcv))&=\big| c_v^{-1}\times L_v(E,1)\big|_p,\\
     \# H^1(F_v,\Ep)&=| c_v|_p^{-1},
    \end{align*}
    which finishes the proof.
   \end{proof}
   
   \begin{remark*}
    \begin{enumerate}[(i)]
     \item Since $J_v(\Fc)$ has finite $\Gamma$-Euler characteristic, it follows from Lemma~\ref{gammamodules2} that the characteristic element of $J_v(\Fc)$ has nonzero constant term of value $| L_v(E,1)|_p$.
     \item By the multiplicativity of characteristic elements in short exact sequences (c.f.~Proposition~\ref{multiplicativity3}), equation~\eqref{largesurjective} implies that
     \[ f_{\SEc}=f_{\Sc}\times\prod_{v\in\fM}f_v, \]
     where for each $v\in\fM$, $f_v$ denotes the characteristic element of $J_v(\Fc)$.
    \end{enumerate}
   \end{remark*}
 
  
\section{Local results}\label{local2}

 For a prime $v$ of $F$, denote by $H_w$ the decomposition group in $H$ of some fixed prime $w$ of $F_\infty$ above $v$; note that $H_w$ can be identified with the local Galois group $\Gal(\Fiw/\Fcw)$. Similarly, let $\Sigma_w=\Gal(\Fiw/F_v)$ and $\Gamma_v=\Gal(\Fcw/F_v)$. We use the notation introduced in Section~\ref{diagram}.
 
 
 \subsection{Local Cohomology}\label{LocCoh2}
  
  \begin{lem}\label{Shapiro2}
   Let $v\in S''$. Then for all $i\geq 0$ there are canonical isomorphisms
   \[
    H^i(H,J_v(F_\infty))\cong \bigoplus_{w| v}H^i(H_w,H^1(F_{\infty,w},E)(p)).
   \]
  \end{lem}
  \begin{proof}
   Easy consequence of Shapiro's lemma.
  \end{proof}
  
  \begin{lem}
   If $v\in S''$ not dividing $p$, then
   \[
    H^i(H,J_v(F_\infty))=0
   \]
   for all $i\geq 0$.
  \end{lem}
  \begin{proof}
   Since the inertia group of $v$ in $\Sigma$ is finite, it follows that $\Fiw$ is a finite extension 
   of $\Fcw$. Now by assumption, $\Sigma_w$ has
   no element of order $p$, so the Galois group
   $H_w=\Gal(\Fiw/\Fcw)$ is of order prime to $p$. The result is now immediate from~\ref{Shapiro2}.
  \end{proof}
  
  Now let $v$ be a prime of $F$ dividing $p$, and denote by $\kiw$ the residue field of $\Fiw$ and by $\tilde{E}_v$ the reduction of $E$ modulo $v$. 
  
  \begin{lem}\label{localisom2}
   Let $v$ be a prime of $F$ dividing $p$ and assume that $E$ has good ordinary reduction at $v$. Then for all $i\geq 1$, there 
   are canonical isomorphisms
   \[ H^i(H_w,H^i(\Fiw,E)(p))\cong H^{i+2}(H_w,\tEp(\kiw)). \] 
  \end{lem} 
  \begin{proof}
   The extension $\Fiw$ of $F_v$ is deeply ramified in the sense of~\cite{coatesgreenberg96} 
   since it contains the deeply ramified field $\Fcw$, so there is a canonical isomorphism of $\Sigma_w$-modules
   \begin{equation}\label{reducisom1}
    H^1(F_{\infty,w},E)(p) \cong H^1(F_{\infty,w},\tE).
   \end{equation}
   Again using Tate local duality, we see that
   \[
    H^i(F_{\infty,w},\tilde{E}_{v,p^\infty})=0
   \]
   for all $i \geq 2$. The lemma now follows from applying Hochschild-Serre for the module on 
   the right of equation~\eqref{reducisom1} 
   to the extension $F_{\infty,w}$ over $\Fcv$.   
  \end{proof}

  
 \subsection{Analysis of the local restriction maps}\label{analoc2}
 
  Recall that we have defined the local restriction map
  \[
   \gamma_S(\Fc)=\bigoplus_{w| S''}\gamma_w,
  \]
  where $w$ runs over all primes of $\Fc$ lying above primes in $S''$, and where
  \[
   \gamma_w:H^1(\Fcw,E)(p) \rightarrow H^1(\Fiw,E)(p)^{H_w}.
  \] 
  
  \begin{lem}\label{localkercokerzero2}
   Let $w$ be a prime of $\Fc$ which divides a prime in $S''$ but does not divide $p$. Then both
   $\ker(\gamma_w)$ and $\coker(\gamma_w)$ are zero.
  \end{lem}
  \begin{proof}
   By inflation-restriction, we have
   \begin{align*}
    \ker(\gamma_w)  =& H^1(H_w,\Ep(\Fiw)),\\
    \coker(\gamma_w)=&H^2(H_w,\Ep(\Fiw)).
   \end{align*}
   But $H_w$ is finite of order prime to $p$, which proves the lemma.
  \end{proof}
  
  \begin{lem}\label{localkercoker2}
   Let $v$ be a prime of $F$ dividing $p$, and let $w$ be a prime of $F_\infty$ 
   above $v$. Then 
   \begin{align*}
    \ker(\gamma_w)  &=H^1(H_w,\tEp(\kiw)),\\
    \coker(\gamma_w)&=H^2(H_w,\tEp(\kiw)).
   \end{align*}
  \end{lem}
  \begin{proof}
   Both $\Fiw$ and $\Fcw$ are deeply ramified in the sense of~\cite{coatesgreenberg96}, so we have 
   isomorphisms
   \begin{align*}
    H^1(\Fcw,E)(p)&\cong H^1(\Fcw,\tEp(\kiw)),\\
    H^2(\Fiw,E)(p)&\cong H^2(\Fiw,\tEp(\kiw)).
   \end{align*}
   The lemma now follows from the inflation-restriction exact sequence. 
  \end{proof}

 
 \section{Global results}\label{global2}
 
  Define the following hypotheses:
  \begin{itemize}
   \item (Tors$_{cyc}$) $\calC (E/\Fc)$ is $\Lambda(\Gamma)$-torsion,
   \item (Tors$_\infty$) $\calC (E/F_\infty)\in\MHG$.
  \end{itemize}

  Recall that (Tors$_\infty$) implies (Tors$_{cyc}$) by Lemma~\ref{Gammatorsion}. To simplify notation, define the condition
  \begin{itemize}
   \item (Fin) Both (Fin$_{glob}$) and (Fin$_{loc}$) are satisfied.
  \end{itemize}
 
 
 
  \subsection{Global cohomology}\label{GlobCoh2}
   
   \begin{lem}\label{globalisom2}
    Suppose that (Tors$_\infty$) holds. Then there are canonical isomorphisms
    \[
     H^i(H,H^1(G_S(F_\infty),\Ep))\cong H^{i+2}(H,\Ep(F_\infty))
    \]
    for all $i\geq 1$.
   \end{lem}
   \begin{proof}
    By Proposition 6 in~\cite{zerbes04}, we have
    \[
     H^i(G_S(F_\infty),\Ep)=0
    \]
    for all $i\geq 2$. Also, it is well-known that the condition that ${\mathcal C}(E/\Fc)$ is 
    $\Lambda(\Gamma)$-torsion implies that
    $H^2(G_S(\Fc),\Ep)=0$. The lemma now follows from applying Hochschild-Serre to the extension 
    $F_\infty$ over $\Fc$.
   \end{proof}
   
   \begin{remark*}
    If $\Ep$  is rational over $F_\infty$, then we do not need assumption (Tors$_\infty$) (c.f.~Theo\-rem 2.10 in~\cite{coates97}).
   \end{remark*}

  
  \subsection{Calculation of $\chi\big(\Gamma,\SE^H\big)$}\label{ECcalculationI2}
  
  In this section, we relate the generalised $\Gamma$-Euler characteristic of $\SE^H$ to the generalised $\Gamma$-Euler characteristic of $\SEc$. Taking  $H$-invariants of~\eqref{Selmer*} gives the commutative diagram

   \begin{diagram}
    0 \to & \SE^H         & \rTo & H^1(G_S(F_\infty),\Ep)^H & \rTo^{\psi_S(F_\infty)} & \bigoplus_{v\in S}J_v(F_\infty)^H &    \\
          & \uTo_{\alpha} &      & \uTo_{\beta}             &                         &  \uTo_{\gamma_S(\Fc)}             &     \\
    0 \to & \SEc          & \rTo & H^1(F_S(\Fc),\Ep)        & \rTo^{\lambda_S(\Fc)}   & \bigoplus_{v\in S''}J_v(\Fc)      & \to 0
   \end{diagram} 
  
   Here, $\gamma_S(\Fc)=\bigoplus_{w| S''}\gamma_w(\Fc)$, where
   \[
    \gamma_w(\Fc): H^1(\Fcw,E)(p)\rTo H^1(\Fiw,E)(p)^{H_w},
   \]
   and we write $w| S''$ if $w$ is a prime of $\Fc$ dividing a prime in $S''$. 
   
  \begin{prop}\label{relation2}
   Let $F_\infty$ be an admissible $p$-adic Lie extension of $F$. Suppose that the conditions (Fin) and (Tors$_{cyc}$) hold. Then $\SE^H$ has finite generalised $\Gamma$-Euler characteristic if and only if so does $\Sel'(E/\Fc)$, and if this is the case, then
   \[ \chi\big(\Gamma,\SE^H\big)=\chi\big(\Gamma,\Sel'(E/\Fc)\big). \]
  \end{prop}
  \begin{proof}
   Apply the snake lemma to the above commutative diagram to get 
   an exact sequence of $\Gamma$-modules
   \begin{align*}
    0&\rTo\ker(\alpha)\rTo\ker(\beta)\rTo\ker(\gamma)\\
     &\rTo\coker(\alpha)\rTo\coker(\beta)\rTo\coker(\gamma)
   \end{align*}
   Let $A=\Sel'(E/\Fc)/\ker(\alpha)$. Inflation-restriction shows that
   \begin{align*}
    \ker(\beta)  &=H^1(H,\Ep(F_\infty)),\\
    \coker(\beta)&=H^2(H,\Ep(F_\infty)),
   \end{align*}
   which are finite by assumption. It follows that $\ker(\alpha)$ 
   is finite, so by Corollary~\ref{sameEC2} $A$ has finite
  generalised $\Gamma$-Euler characteristic if and only if so does $\Sel'(E/\Fc)$, and then
   $\chi(\Gamma,A)=\chi\big(\Gamma,\Sel'(E/\Fc)\big)$. We clearly have an exact sequence of $\Gamma$-modules
   \[
    0\rTo A\rTo\Sel(E/F_\infty)^H \rTo\coker(\alpha)\rTo 0.
   \]
   By Lemma~\ref{sameEC2} it is sufficient to show that $\coker(\alpha)$ has 
   finite $\Gamma$-Euler characteristic, and that
   $\chi\big(\Gamma,\coker(\alpha)\big)=1$. The map $\ker(\gamma)\rTo\coker(\alpha)$ has finite 
   kernel and cokernel, so 
   it is sufficient to show that $\ker(\gamma)$ has finite generalised $\Gamma$-Euler 
   characteristic, and
   $\chi\big(\Gamma,\ker(\gamma)\big)=1$. Now by Lemmas~\ref{localkercokerzero2} and~\ref{localkercoker2}, 
   for $i=0,1$ we have 
   \[
    H^i(\Gamma,\ker(\gamma))\cong\bigoplus_{v| p}H^i(\Gamma_v,H^1(H_w,\tEp(\kiw))),
   \]
   where $w$ is a fixed prime of $\Fc$ above $v$ for each prime $v$ of $F$ dividing $p$. 
   It is therefore sufficient to show that
   $H^1(H_w,\tEp(\kiw))$ has finite generalised $\Gamma_v$-Euler characteristic for each 
   $v| p$, and
   \[
    \chi\big(\Gamma_v,H^1(H_w,\tEp(\kiw))\big)=1.
   \]
   But since we assume that $E$ has good ordinary reduction at all primes of $F$ dividing $p$, 
   $H^i(H_w,\tEp(\kiw))$ is finite for all
   $i\geq 0$ by assumption, so clearly
   \[
    \# H^1(H_w,\tEp(\kiw))^{\Gamma_v} =\# H^1(H_w,\tEp(\kiw))_{\Gamma_v}. 
   \]
  \end{proof}

    
  \subsection{Calculation of $\chi\big(\Sigma,\SE\big)$}\label{ECcalculationII2}
   
   Recall that by Proposition~\ref{surjection}, the assumption (Tors$_\infty$) implies that we have the short exact sequence
   \begin{equation}\label{Selmer2}
    0 \rTo \SE \rTo H^1(G_S(F_\infty),\Ep) \rTo \bigoplus_{v\in S}J_v(F_\infty) \rTo 0.
   \end{equation}
   
  \begin{lem}\label{finitehighercohomology2}
   If the conditions (Tors$_\infty$) and (Fin) hold, then $H^i(H,\SE)$ is finite for all $i\geq 1$.
  \end{lem}
  \begin{proof}
   Recall that we had defined the map
   \[\psi_S(F_\infty):H^1(G_S(F_\infty),\Ep)^H\rTo \bigoplus_{v\in S}J_v(F_\infty)^H.\]
   It is easy to see that $\Im(\gamma)$ is contained in $\Im(\psi_S(F_\infty))$. 
   It is shown in Lemma~\ref{cohomologyatp2} that
   $\coker(\gamma)$ is finite, which implies that $\coker(\psi_S(F_\infty))$ is finite. 
   Taking $H$-cohomology of~\eqref{Selmer2}, we
   get the long exact sequence
   \begin{align*}
    0 &\rTo \coker(\psi_S(F_\infty))\rTo H^1(H,\SE)\rTo H^1(H,H^1(G_S(F_\infty),\Ep))\\
      &\rTo\dots\rTo\bigoplus_{v\in S}H^d(H,J_v(F_\infty))\rTo 0
   \end{align*}
   Here, $d=\dim(H)$ which is equal to the $p$-cohomological dimesion of $H$ since $H$ has no 
   element of order $p$. By
   Lemmas~\ref{globalisom2} and~\ref{localisom2}, for all $i\geq 1$ we have canonical isomorphisms
   \[
    H^i(H,H^1(G_S(F_\infty),\Ep))\cong H^{i+2}(H,\Ep(F_\infty)),
   \]
   and
   \[
    H^i(H,J_v(F_\infty))=
    \begin{cases}
     \bigoplus_{w| v}H^{i+2}(H_w,\tEp(\kiw)) &\text{if $v| p$}\\
     0                                          &\text{otherwise}
    \end{cases}
   \]
   As shown in Lemmas~\ref{globalfinite2} and~\ref{cohomologyatp2}, these cohomology groups 
   are all finite, so we deduce that $H^i(H,H^1(G_S(F_\infty),\Ep))$ and
   $H^i(H,J_v(F_\infty))$ are finite for all $i\geq 1$, which proves the lemma.
  \end{proof}
  
  \begin{prop}\label{finitekercokerd0}
   Let $d_0$ be the map
   \[
    d_0: H^0(\Sigma,\SE)\rTo H^1(\Sigma,\SE)
   \]
   as constructed in Section~\ref{generalizedEC2}. Then under the same conditions as in 
   Lemma~\ref{finitehighercohomology2}, $d_0$ has finite kernel 
   and cokernel if and only if $\Sel'(E/\Fc)$ has finite generalised
   $\Gamma$-Euler characteristic, and then
   \[
    \frac{\#\ker(d_0)}{\#\coker(d_0)}=\frac{\chi\big(\Gamma,\Sel'(E/\Fc)\big)}{\# H^1(H,\SE))^\Gamma}.
   \]
  \end{prop}
  \begin{proof}
   As shown in Proposition~\ref{relation2}, the map
   \[
    f:H^0(\Gamma,\SE^H)\rTo H^1(\Gamma,\SE^H)
   \]
   has finite kernel and cokernel if and only if $\Sel'(E/\Fc)$ has finite generalised 
   $\Gamma$-Euler characteristic, and then  
   \begin{align*}
    \frac{\#\ker(f)}{\#\coker(f)}&=\chi\big(\Gamma,\SE^H\big)\\
                                 &=\chi\big(\Gamma,\Sel'(E/\Fc)\big).
   \end{align*}
   Since $\cd_p(\Gamma)=1$, Hochschild-Serre gives the short exact sequence
   \[
    0\rTo H^1(\Gamma,\SE^H)\rTo H^1(H,\SE)\rTo H^1(H,\SE)^\Gamma\rTo 0.
   \]
   As shown in Lemma~\ref{globalfinite2}, $H^1(H,\SE)$ is finite, which implies the result.
  \end{proof}
   
  \begin{prop}\label{finitesigmaEC2}
   Under the same conditions as in Lemma~\ref{finitehighercohomology2} $\SE$ has finite 
  generalised $\Sigma$-Euler characteristic,
   and
   \[
    \chi\big(\Sigma,\SE\big)=\chi\big(\Gamma,\Sel'(E/\Fc)\big).
   \]
  \end{prop}
  \begin{proof}
   Since $\cd_p(\Gamma)=1$, Hochschild-Serre gives exact sequences
   \[
    0\rTo H^{i-1}(\Gamma,\SE)_\Gamma \rTo H^i(\Sigma,\SE) \rTo H^i(H,\SE)^\Gamma \rTo 0
   \]
   for all $i\geq 2$. It follows from Lemma~\ref{finitehighercohomology2} that 
   $H^i(\Sigma,\SE)$ is finite for all $i\geq 2$. For $i\geq 1$, 
   let
   \[
    h_i=\# H^i(H,\SE)^\Gamma.
   \]
   Since $H^i(H,\SE)$ is finite for all $i\geq 1$, we have
   \[
    \# H^i(H,\SE)^\Gamma =\# H^i(H,\SE)_\Gamma,
   \]
   and so
   \[
    \# H^i(\Sigma,\SE)=h_{i-1}h_i
   \]
   for all $i\geq 2$. Hence
   \[
    \prod_{i\geq 2} (\# H^i(\Sigma,\SE))^{(-1)^i}=h_1.
   \]
   Combining this with Proposition~\ref{finitekercokerd0} shows that $\SE$ has 
   finite generalised $\Sigma$-Euler characteristic, and
   \[
    \chi\big(\Sigma,\SE\big)=\chi\big(\Gamma,\Sel'(E/\Fc)\big).
   \]
  \end{proof}

   
  \subsection{Example}
  
   Let $E=X_1(11)$, which is given by the equation 
   \[
    E:y^2+y=x^3-x^2.
   \]
   Take $p=7$ and $F=\QQ(\mu_7)$. Note that $E$ is ordinary at $7$. Chris Wuthrich has shown the following result.
   
   \begin{prop}
    $\calC(E\slash \Fc)$ is a finitely generated $\ZZ_p$-module of rank $1$. 
   \end{prop}
   
   \begin{cor}
    There exists $u\in\ZZ_p^\times$ such that the characteristic power series of $\calC(E\slash\Fc)$ is $uT$. 
   \end{cor}
   
   \begin{cor}
    $\Sha(E\slash F)(7)$ is finite.
   \end{cor}
   \begin{proof}
    One can show that $\rank_{\ZZ}E(F)=1$. It is well-known that there is a short exact sequence
    \[
     0\rTo \varinjlim E(F)\slash 7^n \rTo \Sel(E\slash F) \rTo\Sha(E\slash F)(7)\rTo 0,
    \]
    so if $\Sha(E\slash F)(7)$ is infinite, then $\rank_{\ZZ_7}\Sel(E\slash F)\geq 2$. Applying the fundamental diagram (c.f. diagram (16) in~\cite{coateshowson97}) to the extension $\Fc$ of $F$, one deduces that the leading term of the characteristic power series of $\calC(E\slash\Fc)$  has exponent $>1$, giving a contradiction.
   \end{proof}

   Let $m>1$ be an integer, and consider the extension $F_\infty=F\big(\mu_{7^\infty},\sqrt[7^\infty]{m}\big)$, which is a Galois extension of $F$ with Galois group $\Sigma\cong\ZZ_p\rtimes\ZZ_p$. 
   
   \begin{lem}
    A prime $l$ of $F$ has infinite inertia group in $F_\infty$ if and only if $l|7m$.
   \end{lem}
   \begin{proof}
    Lemma 3.9 in~\cite{hachimorivenjakob03}.
   \end{proof}
   
  Theorems 3.1 and 3.7 in~\cite{hachimorivenjakob03} imply that  $\calC(E\slash F_\infty)\in\MHG$, so all the conditions of Theorem~\ref{theorem3} are satisfied. Let $m=113$. Since $113=1\mod 7$, $113$ splits completely in $F$. Note that $E$ has good reduction at each prime $v$ of $F$ above $7$ and $113$. By definition, we have
  \[
   L_{v}(E,s) =\big(1+a_vq_v^{-1}+q_v^{-2}\big)^{-1},
  \]
  where $q_v$ is the order of the residue field of $F_v$ and $a_v=q_v+1-\#\tilde{E}_v(k_w)$. Explicit calculation shows that $a_v=-2$ for $v|7$ and $a_v=9$ for all $v|113$, and hence $|L_v(E,1)|_7=7^2$ for $v|7$ and $|L_v(E,1)|_7 = 1$ when $v|113$. It therefore follows from Theorem~\ref{theorem3} that  $\calC(E\slash F_\infty)$ has finite generalised Euler characteristic, and (up to a unit in $\ZZ_7$)
  \[
   \chi\big(\Sigma,\calC(E\slash F_\infty)\big)=7^8.
  \]


\section{Global and local finiteness conditions}\label{examples2}

 \subsection{The general case}
 
  Recall the following conditions:
  
  \begin{itemize}
   \item (Fin$_{glob}$) $H^i(H,\Ep(F_\infty))$ is finite for all $i\geq 0$.
   \item (Fin$_{loc}$) For any prime $v$ of $F$ dividing $p$, $H^i(H_w,\tEp(\kiw))$ is finite for all $i\geq 0$.
  \end{itemize}
    
  The aim of this section is the proof of Theorem~\ref{theorem4}, which defines a class of admissible $p$-adic Lie extensions of $F$ for which the conditions (Fin$_{glob}$) and (Fin$_{loc}$) are satisfied. Assume that $F_\infty$ is the fixed field of the kernel of some continuous representation $(\rho,V)$ of $\Gal(\bar{F}/F)$ on a finite dimensional $\Qp$-vector space $V$. For a prime $v$ of $F$ dividing $p$, denote by $\rho_v$ the restriction of $\rho$ to $\Gal(\bar{F}_v/F_v)$. 
  
  \begin{definition} 
   Let $v|p$, and let $\FF_q$ be the residue field of $F_v$. A Weil number of weight $w\in\ZZ$ with respect to $q$ is an algebraic number all of whose archimedean absolute values are $q^{\frac{w}{2}}$ and its $u$-adic absolute value is $1$ for any non-archimedean prime $u$ which does not divide $p$. 
  \end{definition}

  \begin{prop}\label{cohomologyatp2}
   Let $v$ be a prime of $F$ dividing $p$, and let 
   $w$ be a fixed prime of
   $F_\infty$ above $v$. Assume that the representation $\rho_v$ is semistable. Also, assume that the
   eigenvalues of the Frobenius endomorphism $\Phi$ acting on the filtered module 
   \[D(V)=(\BB_{st}\otimes_{\QQ_p}V)^{\Gal(\bar{F}_v/F_v)}\]
   are Weil numbers of fixed parity. Then $H^i(H_w,\tE(\kiw))$ is finite for all $i\geq 0$. 
  \end{prop}
  \begin{proof}
   If $\tEp(\kiw)$ is finite, then the result is immediate. We can therefore assume that $\tEp$ is rational over $\kiw$.
   Let $H_\infty$ denote the maximal unramified $\ZZ_p$-extension of $\Fcw$ contained in $\Fiw$, and define the Galois group
   $\Delta=\Gal(H_\infty/\Fcw)$. Let $\fH$ and $\fD$ denote the respective Lie algebras. Define the
   $\QQ_p$-vector space $V_p(\tilde{E})=T_p(\tilde{E})\otimes_{\ZZ_p}\QQ_p$. By Theorem 2.4.10 in Chapter V of \cite{lazard65}, it is sufficient to
   show that $H^i(\fH_{\bQp},V_p(\tilde{E})_{\bQp})=0$ for all $i\geq 0$, where for a $\QQ_p$-vector space $W$, we use the 
   notation
   $W_{\bQp}=W\otimes_{\QQ_p}\bQp$. Let $D(V)=(\BB_{st}\otimes_{\QQ_p}V)^{\Gal(\overline{F_v}/F_v)}$, and let $i_1\leq\dots\leq i_d$ denote
   the Hodge-Tate weights of the natural filtration of $D(V)$. Let $\lambda_1,\dots,\lambda_d$ denote the eigenvalues of the
   Frobenius endomorphism $\Phi$ acting on $D(V)$. By assumption, these eigenvalues are Weil numbers of a fixed parity. Then as 
   shown in~\cite{coatessujathawintenberger01}
   there exists an element $X\in\fH$ such that the eigenvalues of $X$ on $V_{\bQp}$ are of the form 
   $\log_\pi(\lambda_1 q^{-i_1}),\dots,\log_\pi(\lambda_d q^{-i_d})$, where the number $\pi$ is chosen to be transcendental over
   $\QQ$. This choice guarantees that $\log_\pi(z)$ is nonzero for every element $z\in\bQp$ which is algebraic over $\QQ$ and 
   not a root of unity. Since the representation of $H_w$ on $V$ is faithful, the eigenvalues $\{\eta_i \}_{i=1}^n$ of
   $\ad_{\fH_{\bQp}}(X)$ are of the form $\log_\pi(\lambda_j q^{-i_j})-\log_\pi(\lambda_k q^{-i_k})$ for some $1\leq j,k\leq d$, so 
   in particular they can be written as $\log_{\pi}(x)$, where $x$ is a 
   Weil number of even weight. By the construction of the element $X$ in~\cite{coatessujathawintenberger01}, it is not difficult to
   see that the image of the projection of $X$ in $\fD_{\bQp}$ is nonzero, so its action on
   $V_p(\tilde{E})_{\bQp}$ is non-trivial. It is well-known that the eigenvalue $\mu$ of this action is of the form $\log_\pi(x)$, where
   $x$ is a Weil number of weight $1$, i.e. it is of absolute complex value $q^{\frac{1}{2}}$. It follows
   that $\mu-(\eta_{j_1}+\dots +\eta_{j_k})$ is of the form $\log_\pi(\kappa)$, where $\kappa$ is a Weil number of complex absolute 
   value $m/2$ with $m$ an odd integer. In particular, we have 
   \[
    \mu-(\eta_{j_1}+\dots +\eta_{j_k})\neq 0
   \]
   for all choices $1\leq j_1<\dots<j_k\leq n$, so the representation of $\fH_{\bQp}$ on $V_p(\tilde{E})_{\bQp}$ satisfies 
   Serre's criterion. As shown in~\cite{coatessujathawintenberger01}, this implies that 
   \[H^i(\fH_{\bQp},V_p(\tilde{E})_{\bQp})=0\] 
   for all $i\geq 0$, which completes the proof.
  \end{proof}
  
  To prove Theorem~\ref{theorem4}, it remains to show that (Fin$_{glob}$) is satisfied. We start by proving the following lemma.
  
  \begin{lem}\label{notrational1}
   Let $F$ be a finite extension of $\QQ$, $E$ an elliptic curve defined over $F$ and $p$ an odd prime. 
   Let $F_\infty$ be a $p$-adic Lie extension of $F$ which contains the cyclotomic 
   $\ZZ_p$-extension. If 
   $\Ep$ is not rational over $F_\infty$, then $\Ep(F_\infty)$ is finite.
  \end{lem}
  \begin{proof}
   Suppose that $\Ep(F_\infty)$ is infinite. Assume first that $E$ does not have complex multiplication. 
   If $\Ep(F_\infty)$ has
   $\ZZ_p$-corank equal to $1$, then $V_p(E)=T_p(E)\otimes_{\ZZ_p}\QQ_p$ has a $1$-dimensional 
   $\Gal(\bar{F}/F)$-invariant subspace. 
   However,  
   as shown in~\cite{serre72}, the Galois group of $F(\Ep)$ over $F$ is an open subgroup of 
   $\GL_2(\ZZ_p)$, so the representation of 
   $\Gal(\bar{F}/F)$ on $V_p(E)$ is irreducible, which gives the required contradiction.
   Assume now that $E$ has complex multiplication. We first show that the $p$-torsion points of $E$ are 
   rational over $F_\infty$. Let 
   $\Delta=\Gal(F_\infty(E_p)/F_\infty)$. It is well-known that the order of $\Delta$ is prime to $p$ 
   (cf.~\cite{perrin-riou84}). By 
   assumption, $E_p(F_\infty)\neq 0$. 
   By choosing a suitable $\ZZ/p\ZZ$-basis of $E_p$, we get an isomorphism of $\Delta$ with a 
   subgroup of $\GL_2(\ZZ/p\ZZ)$ consisting 
   of matrices of the form 
   $\begin{pmatrix}1 &x \\ 0 &1\end{pmatrix}$, where $x\in\ZZ/p\ZZ$. But if $x\neq 0$, then 
   $\begin{pmatrix}1 &x \\ 0 &1\end{pmatrix}$ 
   generates a subgroup of order $p$. It follows that $\Delta$ is trivial and hence $E_p$ is rational 
   over $F_\infty$. It is shown 
   in~\cite{perrin-riou84} that
   $\Gal(F(\Ep)/F)\cong \ZZ_p\times\ZZ_p$. Note that $F(\Ep)$ contains $F^{cyc}$ by the Weil pairing. 
   Let $K_\infty$ denote the extension of $F$ generated 
   by $\Ep(F_\infty)$. Then $\Gal(K_\infty/F)\cong\ZZ_p\times\ZZ/p^n\ZZ$ for some $n\geq 0$. It is shown 
   in~\cite{coatessujatha00} 
   that $\Ep(F^{cyc})$ is finite, so
   $K_\infty$ and $F^{cyc}$ intersect in a finite extension of $F$. It follows that 
   $F(\Ep)=K_\infty F^{cyc}$, which proves the 
   proposition.
  \end{proof}
  
  \begin{prop}\label{globalfinite2}
   Let $F$ be a finite extension of $\QQ$ and $p$ a prime $\geq 5$. Let $F_\infty$ be an admissible
   $p$-adic Lie extension of $F$ with Galois group $\Sigma=\Gal(F_\infty/F)$, and let
   $H=\Gal(F_\infty/\Fc)$. Assume $\Lie H$ is reductive. Let $E$ be an elliptic curve
   defined over $F$. Then $H^i(H,\Ep(F_\infty))$ is finite for all $i\geq 0$. 
  \end{prop}

  \begin{proof}
   If  $\Ep$ is not rational over $F_\infty$, then as shown in Lemma~\ref{notrational1}, $\Ep(F_\infty)$ is finite. Recall the following result from~\cite{serre98-1}, Section 4.1: If $G$ is a pro-$p$ $p$-adic Lie group and $M$ a finite $p$-primary discrete $G$-module, then $H^i(G,M)$ is finite for all $i$. Since $\Ep(F_\infty)$ is finite, this proves the proposition when $H$ is pro-$p$. Suppose that $H$ is not pro-$p$. Fix an open normal pro-$p$ subgroup $G$ of $H$. Since $G$ is pro-$p$, $H^i(G,\Ep(F_\infty))$ is finite for all $i\geq 0$. The finiteness of the $H^i(H,\Ep(F_\infty))$ now follows from the Hochschild-Serre spectral sequence
   \[ H^p(H,H^q(G,\Ep(F_\infty)))\Rightarrow H^{p+q}(H,\Ep(F_\infty)). \]   

   Suppose now that $\Ep$ is rational over $F_\infty$. Let $\fH$ denote the Lie algebra of $H$, and let $V_p(E)=T_p(E)\otimes_{\ZZ_p}\QQ_p$, where $T_p(E)$ denotes the Tate module of $\Ep$. To prove the proposition, it is sufficient to    show that $H^i(\fH,\VE)=0$ for all $i\geq 0$. Let $J=\Gal(F(\Ep)/\Fc)$ and $I=\Gal(F_\infty/\Fc)$. To simply notation, write $\fJ=\Lie J$ and $\fJ'=\Lie I$. Since $\fH$ is reductive, it is easy to see that we have $\fH\cong\fJ'\times\fJ$. Now $H^i(\fJ,\VE)=0$ for all $i\geq 0$  (the case when $E$ has complex multiplication is a direct consequence of Imai's theorem \cite{imai75}, and the non-CM case is shown in the appendix of~\cite{coatessujatha00}). Since in particular $H^0(\fJ,\VE)=0$ and $\fJ'\cong \fH/\fJ$, applying the Hochschild-Serre spectral sequence to the ideal $\fJ$ of $\fH$ finishes the proof.
  \end{proof}



 \subsection{A special case}\label{special2}

  When the $p$-adic Lie extension is generated by the points of $p$-power order on an abelian
  variety, then we can analyse the local cohomology groups in more detail. The aim of this section is to
  show that for such an extension the condition (Fin$_{loc}$) is satisfied, without using the results
  in~\cite{coatessujathawintenberger01}.
  
  Let $E$ be an elliptic curve and $A$
  an abelian variety, both defined over $F$, and assume that both $A$ and $E$ have either good
  ordinary reduction at the the primes of $F$ dividing $p$. Let $F_\infty=F(\Ap)$,
  $\Sigma=\Gal(F_\infty/F)$, $H=\Gal(F_\infty/\Fc)$ and $\Gamma=\Gal(\Fc/F)$.
  
  Let $v$ be a prime of $F$ dividing $p$ and let $w$ be a prime of $\Fc$ above $v$. We split the
  proof of (Fin$_{loc}$) up into two propositions. The first one of these propositions is only important
  here because its proof gives on a detailed study of the local Galois group $H_w$ which we need in order
  to show that (Fin$_{loc}$) is satisfied. 

  \begin{prop}\label{multiplicative3}
   Let $v$ be a prime of $F$ dividing $p$, and let $w$ be a prime of $\Fc$ dividing $v$.  
   Assume that $A$ has good ordinary reduction at $v$. Then 
   \[H^0(H_w,\QZ)\cong H^1(H_w,\QZ)\cong\QZ\] 
   as $\Gamma_v$-modules, and $H^i(H_w,\QQ_p/\ZZ_p)$ is finite for all $i\geq 2$.
  \end{prop}

  \begin{remark*}
   Observe that to prove the finiteness of $H^i(H_w,\QQ_p/\ZZ_p)$ for $i\geq 2$, 
   it is sufficient to show that $H^i(R,\QQ_p/\ZZ_p)$ is finite for some open normal subgroup $R$ of
   $H_w$, i.e. we can replace $F_v$ by any finite Galois extension which is contained in $\Fiw$.
   Let $M_\infty$ be the maximal unramified extension of $\Fcw$ contained in $\Fiw$, and define the 
   Galois groups $M=\Gal(\Fiw/M_\infty)$ and $\Delta=\Gal(M_\infty/\Fcw)$. Then $\Fiw$ is a finite
   extension of the unique unramified $\ZZ_p$-extension of $\Fcw$, so $\Delta$ is isomorphic to the direct
   product of $\ZZ_p$
   with some finite abelian group of order prime to $p$. By the previous remark, we can without loss of 
   generality assume that $\Delta\cong\ZZ_p$.
   Let $\delta$ be topological generator of $\Delta$. As shown in~\cite{greenberg99}, we have
   \begin{align}\label{structureM}
    M&\leq\Hom(T_p(\tilde{A}),T_p(\hat{A}))\\
     &=\Hom(T_p(\tilde{A}),\Hom(T_p(\tilde{A^d}),\ZZ_p))
   \end{align}
   where the equality comes from the Weil pairing and the observation that $H_w$ acts trivially on the
   $p$-power root of unity. Here, $A^d$ denotes the dual abelian variety; note that $A^d$ is isogenous to
   $A$ over $\bar{F}_v$. 
   It is easy to see that the action of $\delta$ on $\Hom(T_p(\tilde{A}),T_p(\hat{A}))$ is
   diagonizable over $\bar{\QQ}_p$, and the eigenvalues are $\{\lambda_i^{-1}\lambda_j^{-1}\}$, where the
   $\lambda_i$ are the eigenvalues of $\Omega$. As shown in~\cite{greenberg99}, the $\lambda_i$ are $p$-adic units and
   Weil numbers of absolute complex value $q^{1/2}$, where $q$ is the cardinality of the residue
   field of $F_v$. 
   Since $M$ is a $\delta$-invariant subspace of $\Hom(T_p(\tilde{A}),T_p(\hat{A}))$, the action of
   $\delta$ on $M$ is also diagonizable. If $\Phi$ denotes the matrix representing this action  
   with respect to some fixed $\ZZ_p$-basis of $M$, then it follows that the eigenvalues 
   $\{\alpha_k\}$ of $\Phi$ are contained in $\{\lambda_i^{-1}\lambda_j^{-1}\}$, so in particular they are $p$-adic units
   and Weil numbers of absolute complex value $q^{-1}$. Now clearly $\Phi\in \GL_n(\ZZ_p)$.
   Since $1+M_n(p\ZZ_p)$ is an open subgroup of $\GL_n(\ZZ_p)$, we may without loss of generality assume
   that $\Phi\in 1+M_n(p\ZZ_p)$ (if necessary, replace $F_v$ by a finite extension contained in $\Fcw$).
   Note that then the eigenvalues of $\Phi$ are elements of $1+\mathfrak{m}$, where $\mathfrak{m}$ is the maximal
   ideal in the ring of integers of $\bar{F_v}$. 
  \end{remark*}
   
  \begin{proof}
   Let $F_v^{nr}$ denote the unramified $\Zp$-extension of $F_v$. It is then clear from the above
   discussion that (up to a finite extension) 
   $\Fcw F_v^{nr}$ is the maximal abelian subextension of $\Fcw$ contained in $\Fiw$, so 
   $H^1(H_w,\QZ)\cong\QZ$, which is an isomorphism of $\Gamma_v$-modules since $\Fcw\cap F_v^{nr}=F_v$. 
   
   Observe that $H^i(H_w,\QQ_p/\ZZ_p)$ is finite if and only if $H^i(H_w,\QQ_p)=0$. Let $\fH$ denote the
   Lie algebra of $H$. Then  $H^i(H_w,\QQ_p)$ is a $\QQ_p$-subspace of $H^i(\fH,\QQ_p)$, so it is 
   sufficient to show that $H^i(\fH,\QQ_p)=0$ for all $i\geq 2$. Denote by $\fD$ and $\fM$ the Lie 
   algebras of $\Delta$ and $M$, respectively. It is easy to see that $\fD=\QQ_p$ and $\fM=\QQ_p^N$, both
   with trivial Lie bracket. Since $M$ is a normal subgroup of $H_w$, it follows that 
   $\fM$ is an ideal of $\fH$, and
   we have a canonical isomorphism $\fD\cong\fH/\fM$. Since $\dim(\fD)=1$, for all 
   $i\geq 2$  Hochschild-Serre~\cite{hochschildserre53} gives  short exact sequences
   \begin{equation}\label{shortexact3*}
    0 \rTo H^1(\fD,H^{i-1}(\fM,\QQ_p))\rTo H^i(\fH,\QQ_p)\rTo H^i(\fM,\QQ_p)^\fD\rTo 0.
   \end{equation}
   Here, the action of $\fD$ on $\fM$ is given by the adjoint representation
   \[
    \ad:\fH\rTo\End(\fM),
   \]
   which factors though $\fD$ since $\fM$ is abelian. The idea now is to relate this action to the action
   by conjugation of $H_w$ on $M$, which factors though $\Delta$ since $M$ is abelian.
   As observed above, $\fM$ is isomorphic to $\QQ_p^N$ with the normal vector addition. If $m_1,\dots,m_N$
   is a $\ZZ_p$-basis of $M$, an element of $M$ is of the form $m_1^{x_1}\dots m_N^{x_N}$ for some
   $x_i\in\ZZ_p$. The $\log$ function is then given by
   \begin{align*}
    \log:                  M&\rTo\fM, \\
    m_1^{x_1}\dots m_N^{x_N}&\rMapsto (x_1,...,x_N).
   \end{align*}
   We remark that this corresponds to the inclusion of groups 
   $M\hookrightarrow\fM$. Again, let $\delta$ be a topological generator of $\Delta$, and denote by
   $c(\delta)$ the action by conjugation on $M$: If $\tilde{\delta}$ is any lift of $\delta$ to $H_w$ and
   $m\in M$, then $c(\delta)(m)=\tilde{\delta}m\tilde{\delta}^{-1}$. If $Lc(\delta)$ denotes the induced
   action on $\fM$, then the following diagram commutes:  
   \begin{diagram}
    c(\delta):  & M      & \rTo & M        \\
                & \dTo   &      & \dTo     \\
    Lc(\delta)	& \fM    & \rTo & \fM
   \end{diagram}       
   Now $Lc(\delta)$ is $\QQ_p$-linear, so if
   we consider our $\ZZ_p$-basis of $M$ as a $\QQ_p$-basis of $\fM$, then $Lc(\delta)\in\Aut(\fM)$ is also
   represented by $\Phi$. The adjoint representation of $\Delta$ is therefore given by
   \begin{align*}
    \Ad:\Delta&\rTo\Aut(\fM),\\
      \delta^x&\rMapsto\Phi^x.
   \end{align*}
   By definition, $\ad=L\Ad$,where
   \[
    L\Ad:\fD\rTo L\Aut(\fM).
   \]
   Since $\fM=\QQ_p^N$, we have $\Aut(\fM)=\GL_N(\QQ_p)$, and so $L\Aut(\fM)=\End(\fM)=M_N(\QQ_p)$ with
   respect to the chosen $\QQ_p$-basis of $\fM$. For a matrix $A=(a_{ij})\in M_N(\QQ_p)$, define 
   \[
    | A|=\sup_{i,j}| a_{ij}|.
   \]
   This norm is compatible with the topology on $M_N(\QQ_p)$, and
   \[
    M_N(p\ZZ_p)=\{A\in M_N(\QQ_p):| A| < p^{-\frac{1}{p-1}}\}.
   \]
   As shown in~\cite{bourbaki73}, we therefore have the exponential function
   \begin{align*}
    \exp:M_N(p\ZZ_p)&\rTo 1+M_N(p\ZZ_p),\\
                   A&\rMapsto \Sigma_{n\geq 0}\frac{A^n}{n!},
   \end{align*}
   which is an analytic isomorphism of $p$-adic analytic manifolds with inverse
   \begin{align*}
    \log:1+M_N(p\ZZ_p)&\rTo M_N(p\ZZ_p),\\
                     A&\rMapsto \Sigma_{n\geq 1}(-1)^{n+1}\frac{(A-1)^n}{n}.
   \end{align*}
   As observed before, we have $\fD=\QQ_p$ with trivial Lie bracket, and the log-function on $\Delta$ is
   given by
   \begin{align*}
    \log:\Delta&\rTo\fD, \\ 
       \delta^x&\rMapsto x.
   \end{align*}
   By the same argument as above, this allows us to identify $\Delta$ with a subgroup of $\fM$, i.e. the
   log-function is the inclusion map $\Delta\hookrightarrow\fD$. We have the commutative diagram 
   \begin{diagram}
    Ad:  & \Delta    & \rTo & 1+M_N(p\ZZ_p)     \\
         & \dTo      &      & \dTo_{\log}       \\
    ad:  & \fD       & \rTo & M_N(p\ZZ_p)
   \end{diagram}       
   Recall that by assumption, we have 
   $\Phi\in 1+M_N(p\ZZ_p)$. Identifying $\Delta$ with its image in $\fD$, it follows that $\ad(\delta)$ is
   represented by the matrix $\log(\Phi)=\Sigma_{i\geq 1}(-1)^{n+1}\frac{(\Phi-1)^n}{n}$. Now $\ad$ is
   $\QQ_p$-linear, so if $x\in\fD=\QQ_p$, then
   \[
    \ad(x)=x\log(\Phi).
   \]
   It follows that the action of $\fD$ on $\fM$ is diagonizable over $\bar{Q_p}$, and the eigenvalues
   $\{\mu_i\}$ of $\ad(\delta)$ are the $p$-adic logarithms of the eigenvalues of $c(\delta)$ acting 
   on $M$. Observe that $\mu_i\neq 0$ for all $i$ since the eigenvalues of $c(\delta)$ are Weil numbers
   and hence cannot be roots of unity. \\
   Let $K$ be a finite extension of $\QQ_p$ over which $\log(\Phi)$ is diagonizable. For a $\QQ_p$-vector
   space $W$, let $W_K=W\otimes_{\QQ_p} K$. It is well-known that
   \[
    H^i(\fH_K,K)=H^i(\fH,\QQ_p)_K
   \]
   for all $i$. It is therefore sufficient to show that $H^i(\fH_K,K)=0$ for all $i\geq 2$. (Recall what
   we are trying to prove!) Tensoring~\eqref{shortexact3*} with $K$ gives short exact sequences
   \[
    0 \rTo H^1(\fD_K,H^{i-1}(\fM_K,K))\rTo H^i(\fH_K,K)\rTo H^i(\fM_K,K)^{\fD_K}\rTo 0
   \]
   for all $i\geq 2$. Here, the action of $\fD_K$ on $\fM_K$ is the natural extension of the action of
   $\fD$  on $\fM$ that we have analysed above, so it is diagonizable. Decompose
   $\fM_K=\fM_1\oplus\dots\oplus\fM_N$ as a direct sum of eigenspaces, and denote by $\mu_i$ the eigenvalue
   of $\delta$ acting on $\fM_i$. Identifying $\fD_K$ with $K$, we have
   shown that an element $k\in \fD_K$ acts on $\fM_i$ as multiplication by $k\mu_i=k\log(\alpha_i)$,
   where $\alpha_i$ is a Weil number of absolute complex value $q^{-1}$. For $0\leq j<N$, let
   $\fH_{>j}=\fH/\fM_{j+1}\oplus\dots\oplus\fM_N$. 
   
   \begin{claim*}[I] Let $X\cong K$ such that (i) $\fM_n$ acts trivially on $X$ for all $1\leq n\leq j$, and
   (ii) $k\in\fD_K$ acts on $X$ as multiplication by $k\log(\alpha)$,where $\alpha$ is a $p$-adic unit and
   Weil number of complex absolute value $q^s$ for some $s>0$. Then $H^i(\fH_j,X)=0$ for all $i\geq 0$.
   \end{claim*}
   
   To prove the claim, we will use repeatedly the following observation: Let $\mathfrak{G}$ be a
   $1$-dimensional Lie algebra over a field $K$, and suppose that $\mathfrak{G}$ acts nontrivially on a
   $1$-dimensional $K$-vector space $V$. Then $V$ has trivial $\mathfrak{G}$-cohomology. 
   
   \begin{proof}[Proof of Claim I] We procede by induction on $j$:\\
   $j=1$: The case when $i=0$ is clear. Suppose that $i>1$. Now $\dim(\fM_1)=1$, so Hochschild-Serre gives
   an exact sequence
   \[
    0\rTo H^1(\fD_K,X)\rTo H^1(\fH_{>1},X)\rTo H^1(\fM_1,X)^{\fD_K}\rTo 0
   \]
   and isomorphisms
   \[
    H^i(\fH_{>1},X)\cong H^{i-1}(\fD_K,H^1(\fM_1,X))
   \]
   for all $i\geq 2$. Now $\dim(\fD_K)=\dim(H^1(\fM_1,X))=1$ and $k\in\fD_K$ acts on $H^1(\fM_1,X)$ as
   multiplication by $k\log(\alpha\alpha_i)$. Now $\alpha\alpha_i$ is a Weil number of complex absolute
   value $q^{1+s}$, so $\log(\alpha\alpha_i)\neq 0$ and hence the $\fD_K$-cohomology of $H^1(\fM_1,X)$ is 
   zero by above observation.\\
   Let $1<j\leq N$ and suppose the claim holds for $j-1$. Now $\dim(\fM_j)=1$, so again Hochschild-Serre
   gives an exact sequence
   \begin{align*}
    0&\rTo H^1(\fH_{>j-1},X)\rTo H^1(\fH_{>j},X)\rTo H^1(\fM_j,X)^{\fH_{>j-1}}\\
     &\rTo\dots\\
     &H^j(\fH_{>j-1},X)\rTo H^j(\fH_{>j},X)\rTo H^{j-1}(\fH_{>j-1},H^1(\fM_j,X))\rTo 0
   \end{align*}
   By hypothesis, $H^i(\fH_{>j-1},X)$ is zero for all $i\geq 0$. 
   Now $\fM_n$ acts trivially on $H^1(\fM_j,X)$ for
   all $1\leq n<j$, and $k\in\fD_K$ acts on it as multiplication by $k\log(\alpha\alpha_j)$. As observed
   before, $\alpha\alpha_j$ is a Weil number of absolute complex weight $q^{1+s}$, so the action is
   nontrivial and hence $H^1(\fM_j,X)$ satisfies conditions (i) and (ii) above. By induction hypothesis, it
   then follows that $H^{j-1}(\fH_{>j-1},H^1(\fM_j,X))$ is zero for all $i\geq 0$, which proves the claim.
   \end{proof}
   
   \begin{claim*}[II] 
    For all $1\leq j<N$, we have $H^i(\fH_j,K)=0$ for all $i\geq 2$.
   \end{claim*}

   \begin{proof}[Proof of Claim II] 
    Again, we procede by induction on $j$.

    $j=1$: We have $\dim(\fH_{>1})=2$, so it is sufficient to prove that $H^2(\fH_{>1},K)=0$.
    Hochschild-Serre gives an exact sequence
    \[
     H^2(\fD_K,K)\rTo H^2(\fH_{>1},K)\rTo H^1(\fD_K,H^1(\fM_1,K))
    \]
    which proves the result using claim I.

    For $1<j\leq N$, the result is proved using Hochschild-Serre and the result of claim I.
    This finishes the proof of claim II and also of the proposition.
   \end{proof}
  \end{proof}
  
  As a corollary, we get the following proposition which gives a new proof of Proposition 2.12
  in~\cite{zerbes04}, using the standard argument as in \text{\it A.2.9} in~\cite{coatessujatha00}.

  \begin{prop}\label{ordinary3}
   Let $v$ be a prime of $F$ diving $p$ where $E$ has good ordinary reduction, and let $w$ be a prime 
   of $\Fc$ dividing $v$. Assume that $A$ has good ordinary reduction at $v$. Then 
   $H^i(H_w,\tEp(\kiw))$ is finite for all $i\geq 0$.
  \end{prop}
  \begin{proof}
   Suppose that $\tEp(\kiw)$ is finite. Let $R$ be an open normal pro-$p$ subgroup of $H_w$.
   It is well-known (c.f. the proof of Proposition 4.2 in~\cite{zerbes04}) that then $H^i(R,\tEp(\kiw))$ is 
   finite for all $i$, and we deduce the
   finiteness of the $H^i(H_w,\tEp(\kiw))$ from the Hochschild-Serre spectral sequence
   \[
    H^p(H_w,H^q(R,\tEp(\kiw)))\Rightarrow H^{p+q}(H_w,\tEp(\kiw)).
   \]
   Suppose now that $\tEp(\kiw)$ is infinite. We will use without comment the notation of the proof of
   the previous proposition. Let $T_p(\tEp)$ denote the Tate module of $\tEp$, and let
   $V=T_p(\tEp)\otimes_{\ZZ_p}\QQ_p$. It is sufficient to prove that $H^i(H_w,V)=0$ for all $i\geq 0$. 
   Now $\delta$ acts on $V$ as multiplication
   by some $p$-adic unit $\nu$, which is a Weil number of absolute complex weight $q^{1/2}$. As before,
   let $K$ be a finite extension of $\QQ_p$ over which $\Phi$ is diagonizable. Then it is sufficient to
   show that $H^i(\fH_K,V_K)=0$ for all $i\geq 0$. The proposition now follows from applying Claim I of
   the proof of Proposition~\ref{multiplicative3} to the $\fH_K$ module $V_K$.
  \end{proof}

\end{document}